\documentclass[a4paper,11pt,final]{article}

\usepackage{amsmath,amsfonts,amssymb,,amsthm,mathptmx,eucal,graphicx}

\usepackage{hyperref}
\hypersetup{unicode = true, extension = pdf, pdftitle = {Sagan Numbers}, pdfauthor = {J. R. G. Mendonça}, pdfsubject={Recreational Mathematics, Number Theory}, pdfkeywords = {Recreational mathematics, science fiction, number theory, normal numbers, BBP-type formul{\ae}}, pagebackref, bookmarks = true, hyperindex = true, hyperfootnotes = true, pdfnewwindow = true, pdftoolbar = true, pdfmenubar = true, pdffitwindow = false, pdfstartview = {FitH}, colorlinks = true, linkcolor = blue, citecolor = blue, urlcolor = blue}

\theoremstyle{definition}
\newtheorem{definition}{Definition}
\newtheorem{exercise}{Exercise}

\def\leq{\leqslant}
\def\geq{\geqslant}

\begin{document}

\pagestyle{empty}
\renewcommand{\thefootnote}{\fnsymbol{footnote}}

\begin{titlepage}

\begin{center}

{\Large \bf Sagan Numbers}

\vspace{6ex}

{\large
J. Ricardo G. Mendon\c{c}a\footnote{The author holds a Ph.D. in Physics (2000) and is currently professor at Universidade Federal de Uberl\^{a}ndia. Email: {\tt \href{mailto:jrgmendonca@hotmail.com}{\nolinkurl{jrgmendonca@hotmail.com}}}.}
}

\vspace{6ex}

{\large \bf Abstract \\}

\vspace{2ex}

\parbox{120mm}
{We define a new class of numbers based on the first occurrence of certain patterns of zeros and ones in the expansion of irracional numbers in a given basis and call them {\it Sagan numbers\/}, since they were first mentioned, in a special case, by the North-american astronomer Carl E. Sagan in his science-fiction novel {\it Contact\/}. Sagan numbers hold connections with a wealth of mathematical ideas. We describe some properties of the newly defined numbers and indicate directions for further amusement. \\

\vspace{2ex}

{\noindent}{\bf Keywords}: Recreational mathematics~$\cdot$ science fiction~$\cdot$ number theory~$\cdot$ normal numbers~$\cdot$ BBP-type formul{\ae}

\vspace{2ex}

{\noindent}{\bf MSC 2000}: 00A08, 97A20, 97A90}

\end{center}

\end{titlepage}


\pagestyle{plain}

\section{Introduction}
\label{intro}

In his novel {\it Contact\/} \cite{contact}, the North-american astronomer and popularizer of scince Carl E. Sagan (1934-1996) tells the story of a very bright astrophysicist, Dr. Eleanor Ann Arroway, that in her looking for extraterrestrial intelligence, a search that Sagan himself pursued intensely, identifies a message from Vega, the alpha star from the constellation of Lyr{\ae} some 25.3 light-years ($\sim$239.4 trillion kilometers) away from our Sun. After decoded, the message instructs for the building of a machine of unknown purpose that, later, is identified as a five-seats transport supposedly designed to take its occupants to Vega. After some misfortune, eventually the machine is built and Dr.~Arroway, together with four other fellows, finally travels the Milky Way through various space-time highways (warmholes and possibly other types of exotic space-time structures) to meet the senders of the message.

Expectadly, in the encounter various questions as to what the life, the universe and everything mean were addressed. In a remarkable dialog, Dr.~Arroway asks the being with which she met what their experiences with the numinous are like. Quoting the book:

\begin{quote}{\it
``I want to know about your myths, your religions. What fills you with awe? Or are those who make the numinous unable to feel it?''

``(\ldots) Certainly we feel it. (\ldots) I don't say this is it exactly, but it'll give you a flavor of our numinous. It concerns pi, the ratio of the circumference of a circle to its diameter. (\ldots) Our mathematicians have made an effort to calculate it out to\ldots\ none of you seem to know\ldots\ Let's say the ten-billionth place. You won't be surprised to hear that other mathematicians have gone further. Well, eventually---let's say it's in the ten-to-the-twentieth-power place---something happens. The randomly varying digits disappear, and for an unbelievably long time there's nothing but ones and zeros.''

``And the zeros and ones finally stop? You get back to a random sequence of digits? And the number of zeros and ones? Is it a product of prime numbers?''

``Yes, eleven of them.''

``You're telling me there's a message in eleven dimensions hidden deep inside the number pi? Someone in the universe communicates by\ldots math\-e\-mat\-ics? (\ldots) Mathematics isn't arbitrary. I mean pi has to have the same value everywhere. How can you hide a message inside pi? It's built into the fabric of the universe.''

``Exactly.''
}\end{quote}

Upon her return to Earth, and after some more misfortunes but also some glory, Dr.~Arroway ends up setting a program for the searching of extraterrestrial intelligence. One of the projects of the program consisted in assembling huge computational powers to compute $\pi$ in various bases until some anomaly shows up.

One day, the computers start to bip frantically: after reaching more than $10^{20}$ digits after the decimal point in the expansion of $\pi$ in base 11---which is, by no coincidence, the dimensionality of the space-time fabric according to modern superstring theories \cite{greene,smolin}---, the programs figure out a special pattern of 0's and 1's out of the otherwise apparently random digits that, when arranged as a square of specific dimensions, forms a circle. In the words of Carl Sagan, 

\begin{quote}{\it
The program reassembled the digits into a square raster, an equal number across and down. The first line was an uninterrupted file of zeros, left to right. The second line showed a single numeral one, exactly in the middle, with zeros to the borders, left and right. After a few more lines, an unmistakable arc had formed, composed of ones. The simple geometrical figure had been quickly constructed, line by line, self-reflexive, rich with promise. The last line of the figure emerged, all zeros except for a single centered one. The subsequent line would be zeros only, part of the frame.

Hiding in the alternating patterns of digits, deep inside the transcendental number, was a perfect circle, its form traced out by unities in a field of noughts.
}\end{quote}

In the novel, the length of the string of 0's and 1's is said to be given by the product of 11 prime numbers, without stating which prime numbers. Unfortunately, there can not be an integer square root of a product of any eleven primes $p_i$ such that a string of $p_{1} \cdot p_{2} \cdot p_{3} \cdots p_{10} \cdot p_{11}$ symbols can be ``reassembled (\ldots) into a square raster, an equal number across and down.'' So, there is a certain {\it licentia poetica} in Sagan's text. Anyway, supposing that the prime numbers are the first eleven different ones, we would have a string of size $2 \cdot 3 \cdot 5 \cdots 29 \cdot 31 = 200,560,490,130 \simeq 447,840^2$ bits, a very long string indeed. Otherwise, supposing they are just the first prime number repeated eleven times, we would have the much more modest number $2 \cdot 2 \cdot 2 \cdots 2 \cdot 2 = 2^{11} = 2048 \simeq 45^2$, although in this case the sequence of ``nothing but ones and zeros'' could hardly be described as during ``an unbelievably long time.''

In this article we use Sagan's ``$\pi$ revelation'' in {\it Contact\/} to define a new class of integer numbers based of the identification of certain patterns in the $b$-ary expansions of irrational numbers. We call these numbers {\it Sagan numbers\/}. Since most $b$-ary expansions of most irrational numbers are infinite, Sagan numbers form an infinite family of numbers. This article goes as follows: in section~\ref{numbers} we define the Sagan numbers through an ancillary entity called digital $n$-circles and briefly discuss the concept of normal number and what it has to do with Sagan numbers; in section~\ref{directions} some directions for further amusement with the newly defined numbers and their relationship with other branches of mathematics are indicated; and in section~\ref{conclusion} we conclude the article with some general remarks.

There are some exercises scattered throughout the article. Eleven of them! Some are pretty elementary, some are more challenging. We believe that all exercises can be tackled by an average undergraduate student in any scientific discipline. The reader is invited to spend a little time on them to get a better understanding of some of the mathematics related with the subject. As usual, mathematics is learnt by doing, not by watching other people do.


\section{Digital $n$-circles and Sagan numbers}
\label{numbers}

Firstly we need to define certain sequences of 0's and 1's that will embody the patterns of interest to the study of Sagan numbers. We then proceed with the main objects of this article and explore some of their properties.

\subsection{Digital $n$-circles\label{n-circles}}

As one could easily guess from the excerpt from Sagan's novel in the introduction, we are interested in defining patterns of 0's and 1's that when arranged into a square array resemble circles. We call these patterns {\it digital $n$-circles\/}.

Although digital $n$-circles are ancillary in our context, they have many interesting properties and have already been studied (maybe under other names, mostly in unnamed form) in a number of different fields, e.g., in number theory, algebraic geometry, combinatorics, and, of course, computer graphics \cite{graphics}.

\begin{definition}[Digital $n$-circle]
A digital $n$-circle $\Gamma^{(n)} = \langle \gamma_1, \gamma_2, \ldots, \gamma_{n^2} \rangle$ is a string of $n^2$ bits $\gamma_i = 0$ or 1 that when assembled from left to right and from top to botton into an $n \times n$ square raster encodes the pattern of a digitized circle of diameter $n$.
\end{definition}

Compared with the ``definition'' given in {\it Contact\/}, our definition of $\Gamma^{(n)}$ ignores the frame of zeros around the digitized circle. If in the one hand this frame is immaterial for our purposes, on the other hand if each digital $n$-circle is to be accompanied by such a frame, the pattern gains $4n-4$ more bits and overall Sagan numbers become even more rarified than they already are, or, equivalently, become much larger integers; {\it cf.\/} the discussion following.

There are many ways to digitize a circle over a square raster of $n \times n$ pixels to obtain its representation in terms of a string of $n^2$ bits or, alternatively, as $n$ binary vectors of lenght $n$. A first approach would be to draw the curve inside the square raster and to attribute a bit 1 to every square that is crossed by a segment of the circle and a 0 to the noncrossed squares. For circles, it gives the following first digital $n$-circles: $\Gamma^{(1)} = \langle 1 \rangle$, $\Gamma^{(2)} = \langle 1,1;\, 1,1 \rangle$, and $\Gamma^{(3)} = \langle 1,1,1;\, 1,0,1;\, 1,1,1 \rangle$, where for the sake of readability we added a semicolon and a little space between diferent groups of digits corresponding to different lines of the raster bitmap. However, for circles this approach is somewhat ``overshooting,'' crossing more pixels than it is strictly necessary. A possible alternative is to use a circle of diameter $n-1$ instead of diameter $n$. This gives a ``less square'' representation of the circle for small $n$, but does not make much difference for larger diameters. This possibility, however, has an interest in its own not only in dimension two but also in higher dimensional spaces, and is related with a number of geometric and combinatorial questions that arise in fields as diverse as number theory and statistical mechanics.

\begin{exercise}
\label{ex-digitize}
(a)~Verify that in order to digitize a circle {\it in the scheme described above\/} such that the first line shows ``a single numeral one, exactly in the middle, with zeros to the borders, left and right'' we would have to have an odd $n < 5/4$, that is, an impossibility except in the trivial case $n=1$. (b)~Verify that a circle of diameter $n$ does not cross the corner pixel of a square raster of $n \times n$ pixels for $n > 4+2\sqrt{2} \simeq 6.828$. Therefore, the ``square'' aspect of $\Gamma^{(n)}$ smoothes progressively as $n$ gets larger than 7. (c) Show that to obtain ``a single numeral one, exactly in the middle, with zeros to the borders, left and right'' in the first line of the raster bitmap the circle should have a radius $ n/2-1 < r < \sqrt{(n/2-1)^2+(1/2)^2}$, with $n \geq 3$ odd.
\end{exercise}

\begin{exercise}[{\cite[Exercise 3.33]{concrete}}]
\label{ex-concrete}
A circle of diameter $2n-1$ is drawn symmetrically on a $2n \times 2n$ chessboard. (a)~How many cells of the board contain a segment of the circle? (b)~Find a function $f(n,k)$ such that exactly $\sum_{k=1}^{n-1}f(n,k)$ cells of the board lies entirely within the circle. {\it Hint\/}: Notice that $r^2$, with $r$ the radius of the circle, cannot be an integer, so by Pythagoras' theorem the circle does not pass through the corner of any cell.
\end{exercise}

\begin{figure}[tp]
\centering
\includegraphics[viewport= 124 51 521 792, scale=0.45, angle=-90]%
{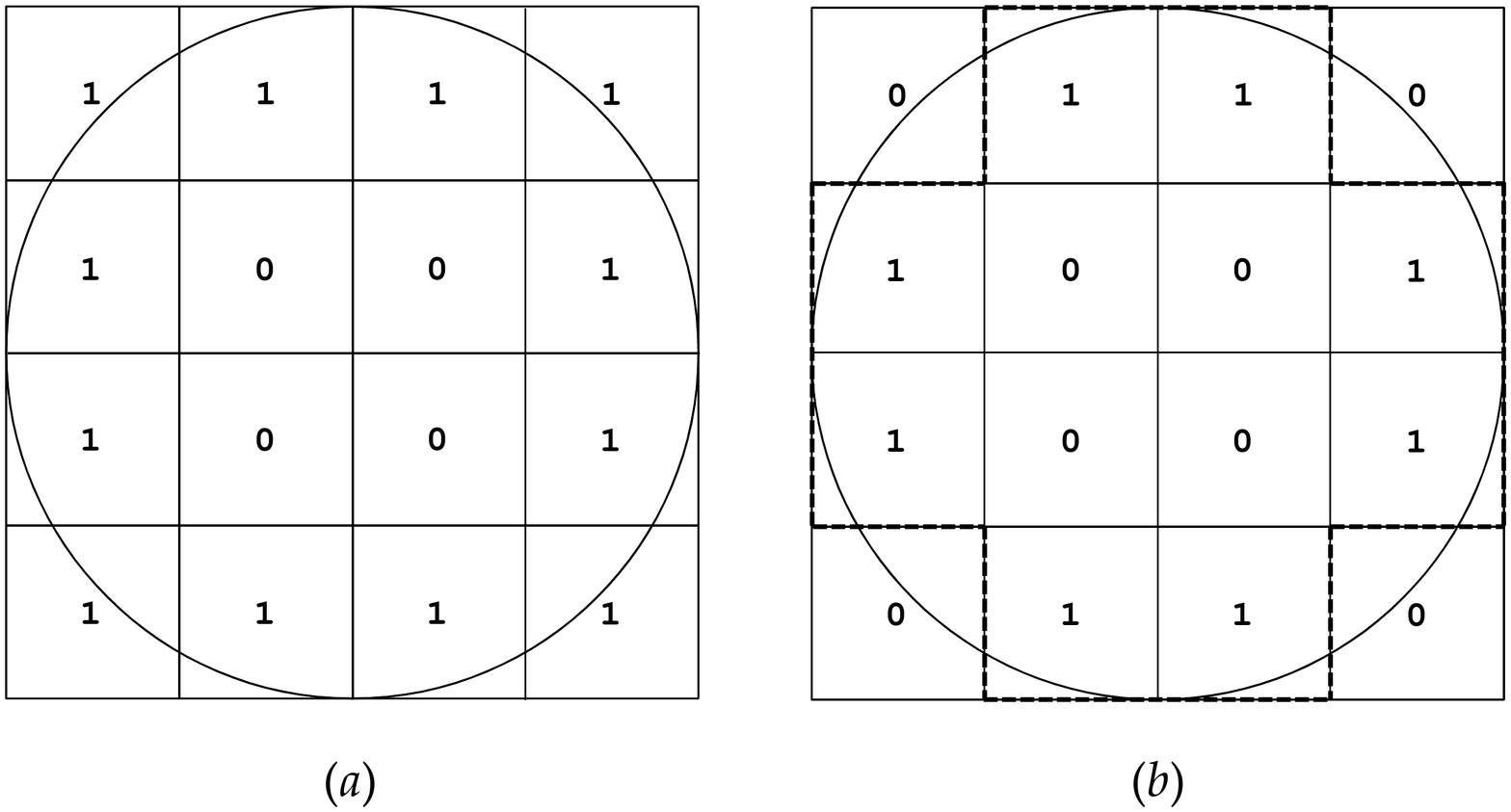}
\caption{Different ways to digitize a circle as a $4 \times 4$ square raster. ($a$) ``Na\"{\i}ve'' approach, giving too much a ``square'' digitized circle for small $n$; ($b$) The preferred approach in raster graphics technology and algorithms.}
\label{fig1}
\end{figure}

The preferred representation in computer graphics to digitize conic sections (circles, ellipses, parabolas, and hyperbolas, together with the degenerate conic, a straight line) is to focus on the edges between adjacent pixels instead of on the pixels themselves. Given a paramatrized curve $z(t) = (x(t), y(t))$, $t \in [0,1] \subset \mathbb{R}$, we digitize it by forcing $z(t)$ to take discrete steps along the edges of the  pixels and by defining which pixels lie inside or outside of the curve. The first pixels that, according to the given definition (which may vary), lie inside of the curve from the outside in are marked with a 1, the other pixels receive a 0. For a circle of radius $r$, $z(t) = (r\cos 2\pi t, r\sin 2\pi t)$ and we can eliminate $t$ to obtain the non-parametric equation $x^{2} + y^{2} = r^{2}$. We then define an inside pixel if its center coordinates, given by $(x-\frac{1}{2}, y-\frac{1}{2})$, observe $(x-\frac{1}{2})^{2} + (y-\frac{1}{2})^{2} \leq r^{2}$.  There are a number of well known, efficient algorithms to implement this approach \cite{pitteway,vanaken,knuth4}. Figure \ref{fig1} illustrates two options on digitizing a circle in a $4 \times 4$ square array, one more ``na\"{\i}ve,'' marking all pixels crossed by the circle of diameter $n$, and the preferred option in raster graphics technology that focus on the edges of pixels whose center lies inside the circle.

\begin{exercise}
\label{ex-algorithm}
Digitize the ellipse $(x/4)^2 + (y/3)^2 = 1$ by hand focusing on the edges and centers of the pixels as described in the text. {\it Hint}: Look up Algorithm T in \cite[Sec.~7.1.3]{knuth4}.
\end{exercise}

\begin{figure}[tp]
\centering
\includegraphics[viewport= 44 42 417 797, scale=0.45, angle=-90]%
{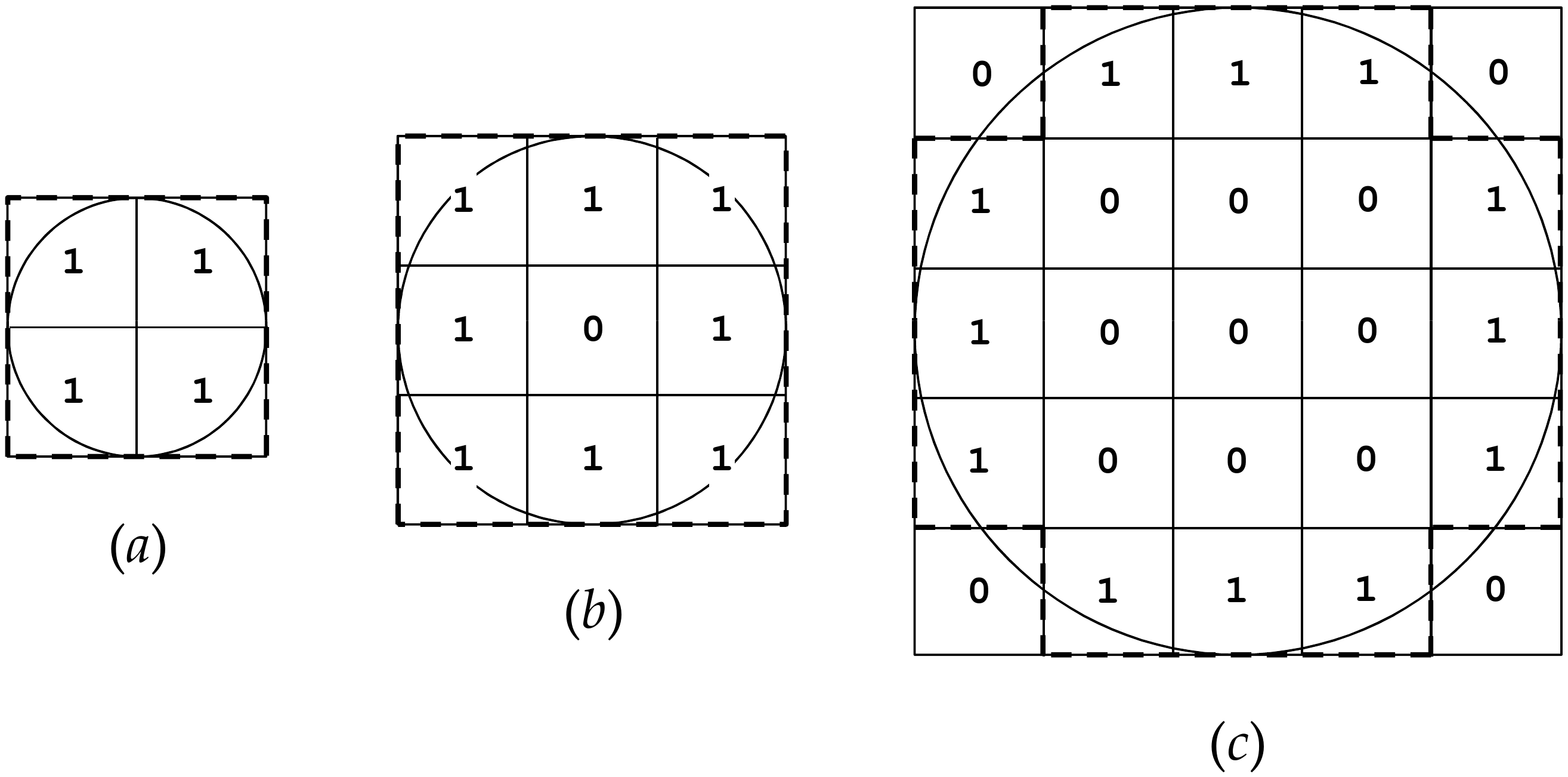}
\caption{Defining patterns for the digital $n$-circles ($a$) $\Gamma^{(2)}$, ($b$) $\Gamma^{(3)}$, and ($c$) $\Gamma^{(5)}$. The defining pattern of $\Gamma^{(4)}$ appears in Fig.~\ref{fig1}($b$).}
\label{fig2}
\end{figure}

\subsection{Sagan numbers}
\label{sagan}

\begin{definition}[Sagan number]
The $n$th decimal Sagan number $S_{\;11}^{(n)}(\pi)$ is the position in the fractional part of the 11-ary expansion of $\pi$ at which the first digit of the digital $n$-circle pattern $\Gamma^{(n)}$ anchors.
\end{definition}

Sagan numbers $S_{\;11}^{(n)}(\pi)$ can be generalized in two obvious ways, as its very notation suggests. The first one is to allow general irrational numbers $\alpha$ in the place of $\pi$, and the other one is to look for digital $n$-circles in arbitrary basis $b$ other then 11. In either way, the generalized definition looses its connection with Sagan's {\it Contact}. The generalized Sagan numbers, however, present very interesting possibilities for recreational mathematics. We thus define the {\it generalized Sagan number\/} as follows:

\begin{definition}[Generalized Sagan number]
The $n$th decimal generalized Sagan number $S_{\;b}^{(n)}(\alpha)$ is the position in the fractional part of the $b$-ary expansion of $\alpha \in \mathbb{R}$ at which the first digit of the digital $n$-circle pattern $\Gamma^{(n)}$ anchors.
\end{definition}

Sagan numbers are, except for the first few ones, very large integers. We can actually estimate how large Sagan numbers $S_{\;11}^{(n)}(\pi)$ should be. If all patterns of $n^2$ digits in the expansion of $\pi$ in base 11 were equally probable (which would amount at saying that $\pi$ is normal to base 11, see next subsection) then we should expect to find, with a finite probability, any particular pattern, say a digital $n$-circle, once every $11^{n^2}$ digits are examined. For a base $b$ other than 11, just change $b$ for 11: we should expect to find any pattern of $n^2$ digits with a finite probability once every $b^{n^2}$ digits are looked up. We then see that while we would most probably find $S_{\;11}^{(1)}(\pi)$ is less than, say, 10 or 20, $S_{\;11}^{(2)}(\pi)$ would most probably be in the range of 10,000--20,000 and $S_{\;11}^{(3)}(\pi)$ in the range of 2,000,000,000--10,000,000,000. Indeed, in base 10 we find that $S_{\;10}^{(1)}(\pi)=1$, the `1' just after the trailing integer `3' in 3.{\underline 1}41\,592\ldots, and $S_{\;10}^{(2)}(\pi) = 12,700$, in a fragment that reads \ldots 144\,{\underline{111\,1}}26\ldots\ In base 11, $S_{\;11}^{(1)}(\pi) = 1$ too, since $\pi_{11}$ = 3.{\underline 1}61\,507\ldots Notice, however, that this estimation is statistical in nature. The first 0 in the decimal expansion of $\pi$ occurs only in position 32!

In Sagan's novel, the pattern possesses {\it at least\/} $2^{11} = 2048$ digits, the smallest product of eleven prime numbers possible, and this means that we would have to examine $11^{2048} \simeq 5.919 \times 10^{2132}$digits before we could possibly find the corresponding $\Gamma^{(n)}$, in this case with $n = \sqrt{2048} \simeq 45$. Since there are less than $3.2 \times 10^{16}$ nanoseconds in a year (even in a leap year), it would take ${\sim}10^{2106}$ times the age of our Universe, currently estimated at about 13.5 billion years, of CPU time to perform the search for $S_{\;11}^{(45)}(\pi)$. Unless there are some radically, otherworldly different way of doing mathematics or computation, this is not going to happen. Actually, there are physical reasons beyond the simple time accounting to believe that this is not going to happen at all, whichever technology or approach, human or alien, one adopts toward numbers and computation in the known Universe \cite{lloyd}.

\begin{exercise}
\label{ex-huge}
Supposing $\pi_{11}$ normal, estimate the amount of CPU time (considering $\sim$1 flop/ns, as above) it would most probably take to find {$S_{\;11}^{(447,840)}(\pi)$}, where $447,840 \simeq \sqrt{2 \cdot 3 \cdot 5 \cdots 29 \cdot 31}$, the number that Carl Sagan probably had in mind when he wrote {\it Contact\/}.
\end{exercise}

\begin{exercise}
\label{ex-base11}
Determine $S_{\;11}^{(2)}(\pi)$. Useful sources of digits of $\pi$ in base 10 enough for the conversion to base 11 can be found on the Internet---just google ``pi digits'' for them. Be careful with the input in your algorithm to avoid the infamous garbage-in-garbage-out ``phenomenon.''
\end{exercise}

\subsection{Normal numbers}
\label{normal}

A real number $\alpha \in \mathbb{R}$ is normal in base $b \in \mathbb{N}$, $b \geq 2$ if in its representation in base $b$ each of the $b^{k}$ different strings $s_k \in \{0, 1, \ldots, b-1\}^k$ occur, in an asymptotic sense, equally often. In other words, $\lim_{l \to \infty}N(s_k,l)/l = b^{-k}$ for each $s_k$, where $N(s_k,l)$ is the number of occurrences of $s_k$ in the first $l$ digits of $\alpha$ in base $b$. A number that is normal in all (integer) bases is called absolutely normal. Obviously, only irrational numbers $\alpha \in \mathbb{R} \backslash \mathbb{Q}$ can be normal.

The normality of an irrational number is a tough question, and few irrational numbers have been proved normal to a given base to date. A very important result in this respect is the one by \'{E}mile Borel that, in 1909, besides introducing the very concept of normal numbers, proved that as a consequence of the strong law of large numbers in probability theory almost every real number is normal to any basis \cite{borel}. The proof is not constructive, however, and it remains an elusive matter to prove a given number normal. On the constructive side, two widely given examples of decimal normal numbers are Champernowne's constant $0.1\,2\,3\,4\,5\,6\,7\,8\,9\,10\,11\ldots$, obtained by concatenating the natural numbers in base 10 in order \cite{champer}, and the Copeland-Erd\H{o}s con\-stant $0.2\,3\,5\,7\,11\,13\,17\ldots$, obtained by concatenating the prime numbers in base 10 in order \cite{erdos}. The binary Champernowne constant $0.1\,10\,11\,100\,101\ldots$ is also normal to base 2 \cite{crandall}. Similar artificial numbers can be constructed in other bases \cite{stone,pincus}, while known irrational numbers like $\pi$ or $\zeta(3)$, the Riemann zeta function at $z=3$, could not yet be proven normal to any basis. Despite the fact that it does appear, based on recent very large calculations that have computed more than $2.576 \times 10^{12}$ decimal digits (a little bit less than $75 \times 2^{35}$ decimal digits, after some discarding) of both $\pi$ and $1/\pi$ that $\pi$ is normal to bases 10 and 16 \cite{kanada,daisuke}, it remains an open question whether $\pi$ is absolutely normal \cite{wagon,quest,calogero}. The quest for normal numbers is an active field in higher arithmetic and number theory, and, more modernly, in what has been dubbed ``experimental mathematics'' \cite{crandall,quest,misiu,kapoor}.

If $\pi$ is absolutely normal, Sagan's argument in {\it Contact}---that the cosmic engineers creators of this Universe have signed their creation by inserting a recognizable pattern in one of its fundamental constants---would be somewhat smeared, for in this case we should expect to find {\it any\/} given Sagan number $S_{\;11}^{(n)}(\pi)$ given enough time or computational power. Now, if $\pi$ is found to be normal only to base 11 or, otherwise, except to base 11, that would be thrilling! We mention in passing that the idea of finding a signature of the creator has also been explored by Douglas Adams in his books {\it The Hitchiker's Guide to the Galaxy\/} and {\it The Restaurant at the End of the Universe} \cite{adams}, in which the ``venerable Magrathean planetary designer'' character Slartibartfast leaves his signature somewhere in Earth's Norwegian {\it fjords\/}, an awarded design piece of him.\footnote{Douglas Adams once said in an interview that he started with ``Phartiphukborlz'' for the name of this character and changed it gradually until he had ``something which sounded that rude, but was almost, but not quite, entirely inoffensive'' {\it and\/} could be broadcast by the BBC.}

\begin{exercise}
\label{ex-random}
In a random sequence of 1,000 digits, what is the probability that there are {\it exactly\/} 100 of each possible digit 0, 1, \ldots, 9? What is the probability that there are {\it exactly\/} 10 of each possible two-digit number 00, 01, \ldots, 99? {\it Hint\/}: Use Stirling's approximation $n! \approx \sqrt{2\pi n}\,(n/e)^n$ to evaluate the large multinomial coefficients that will appear.
\end{exercise}

\begin{exercise}
\label{ex-fibonacci}
Let $\mathcal{F} = 0.0\,1\,1\,2\,3\,5\ldots = 0.F_{0}\,F_{1}\,F_{2}\,F_{3}\,F_{4}\,F_{5}\ldots$, with $F_{0}=0$, $F_{1}=1$, and $F_{k}=F_{k-1}+F_{k-2}$ for $k \geq 2$ the Fibonacci numbers. Prove or disprove that $\mathcal{F}$ is normal to base 10. {\it Hint:\/} Copeland and Erd\H{o}s may give you ideas \cite{erdos}.
\end{exercise}

\begin{exercise}
\label{ex-cfraction}
Instead of concatenating the Fibonacci numbers as in the previous exercise, form the continued fraction
\[
\mathcal{F} = F_0 + \cfrac{1}{F_1+\cfrac{1}{F_2+\cfrac{1}{F_3+\ddots}}} = [F_0; F_1, F_2, F_3, \ldots].
\]
What can you say about the normality of $\mathcal{F}$? Try some computer experiments to obtain the first few thousand decimal places of $\mathcal{F}$ and analyse the relative frequencies of strings of up to 2 or 3 digits as a crude test for normality.  Evaluate the statistical significance of your results against the hypothesis that $\mathcal{F}$ is 10-normal.
\end{exercise}


\section{Directions for further amusement}
\label{directions}

The definitions of $\Gamma^{(n)}$ and $S_{\;b}^{(n)}(\alpha)$ poses some interesting questions to the recreational mathematician. Here we list only a couple of them.

\subsection{Point groups in two dimensions}

Is it possible to obtain the string of 0's and 1's of the digital $n$-circle given only $n$? This may be viewed as an interesting question in analytic geometry. Lots of symmetries exist in the digital $n$-circle string. The string of the digital $n$-circle is also a palindrome (also known as ``capicua''), i.e., it is symmetric about its midpoint and reads the same forward or backward. For $n$ even, $\Gamma^{(n)}$ is of the form 
\[
\Gamma^{(n)} = \Gamma_{\;1}^{(n)} \vee \Gamma_{\;2}^{(n)} \vee \ldots \vee \Gamma_{\;\frac{n}{2}}^{(n)} \vee \Gamma_{\;\frac{n}{2}}^{(n)} \vee \ldots \vee \Gamma_{\;2}^{(n)} \vee \Gamma_{\;1}^{(n)},
\]
where each $\Gamma_{\;k}^{(n)}$ is the pattern given by the bits $\langle \gamma_{(k-1)n+1}, \gamma_{(k-1)n+2}, \ldots, \gamma_{kn} \rangle$ and $\vee$ concatenates strings. This means that it suffices to know the upper semicircle to determine the circle completely. Now, each $\Gamma_{\;k}^{(n)}$ observes a specular (mirror) symmetry about its midpoint, $\Gamma_{\;k}^{(n)}(l) = \Gamma_{\;k}^{(n)}(n-l+1)$. Viewed as a square $n \times n$ matrix, it is also clear that $\Gamma^{(n)}$ is symmetric for all $n$, even or odd, i.e., that $\Gamma_{\;k}^{(n)}(l) = \Gamma_{\;l}^{(n)}(k)$.  This, together with the previous symmetry, means that it suffices to know one quadrant of the circle to determine the whole pattern, the remaining quadrants being easily covered by symmetry relations. Up to the segmentation in quadrants we benefit from the square symmetry of the raster bitmap, but the symmetries become more complicated when it comes to dividing $\Gamma^{(n)}$ in more parts. However, we can push this game further and successively conclude that only one half-quadrant ($0 \leq \theta \leq \pi/4$) is needed to determine the circle completely, then that only a quarter-quadrant ($0 \leq \theta \leq \pi/8$) is needed, and so on. Indeed, to fully determine a circle, we only need two points: its center and a point anywhere on the circle to fix the radius. However, because of the finiteness of $n$ it is clear that the thinning process cannot be carried out indefinitely for digital $n$-circles.

\begin{exercise}
\label{ex-octant}
Determine the relationships between $\Gamma_{\;k_1}^{(n)}(l)$ in the first half-quad\-rant with the other $\Gamma_{\;k_i}^{(n)}(l)$, $2 \leq i \leq 8$, in the other half-quadrants. What is the smallest $n$ below which this cannot be done? {\it Hint}: You may want to draw a few $\Gamma^{(n)}$ in squared paper before proceeding with this exercise.
\end{exercise}

In summary, the subject provides many opportunities to play with the notion of point groups in two dimensions. The implementation of these symmetries in a program to rasterize geometric figures may prove instructive and useful for the student of analytic geometry, computer graphics, and related matters.

\begin{exercise}
\label{ex-program}
Write a computer program that given an integer $n$ outputs $\Gamma^{(n)}$. Implement as much of the symmetry of $\Gamma^{(n)}$ as you can in your program.
\end{exercise}

\subsection{Generalized digital $n$-circles}

In this article, we proposed both the Sagan numbers $S_{\;11}^{(n)}(\pi)$ and their generalization $S_{\;b}^{(n)}(\alpha)$ to other irrational arguments $\alpha$ and bases $b$. It is possible to generalize the digital $n$-circles, too. Sagan described the $\Gamma^{(n)}$ inside of $\pi$ as ``{\ldots}a perfect circle, its form traced out by unities in a field of noughts.'' We could relax both the parts ``traced out by unities'' and the part ``in a field of noughts'' in the definition of $\Gamma^{(n)}$ by letting the individual bits (now, digits) $\gamma_i$ in $\Gamma^{(n)}$ assume arbitrary values $p$ and $q$, such that $\Gamma^{(n)}$ describes a circle traced out by $p$'s in a field of $q$'s. Ultimately, we can define $p$ and $q$ to be nonempty sets $P = \{p_1, p_2, \ldots, p_k\}$, the set of digits that can appear in the rasterized circle, and $Q = \{q_1, q_2, \ldots, q_l\}$, the set of digits that can form the background field. Notice that, although desirable, it is not strictly necessary that $P \cap Q = \varnothing$. We then define the generalized digital $n$-circle as follows:

\begin{definition}[Generalized digital $n$-circle]
A generalized digital $n$-circle {$\Gamma^{(n)}$ $(P,Q)$} = $\langle \gamma_1, \gamma_2, \ldots, \gamma_{n^2} \rangle$ is a string of $n^2$ digits $\gamma_i \in P \cup Q$ that when assembled from left to right and from top to botton into an $n \times n$ square raster encodes the pattern of a digitized circle of diameter $n$ traced out by $\gamma_i \in P$ in a field of $\gamma_i \in Q$.
\end{definition}

The digital $n$-circle $\Gamma^{(n)}$ corresponds to the special case $\Gamma^{(n)}(\{1\},\{0\})$, which we may called the ``plain'' or the ``minimal'' digital $n$-circle. The other extreme is given by the generalized digital $n$-circle $\Gamma^{(n)}(P,Q)$ with $P \subset Q$, that hardly conveys any information and we can call the ``degenerate digital $n$-circle.'' The case $\Gamma^{(n)}(P,Q)$ with $Q \subset P$ is also highly undesirable, although since there are much more background digits $\gamma_i \in Q$ than circle digits $\gamma_i \in P$ in each pattern, this case still conveys some information on the rasterized circle. It is possible to precise the meaning of ``conveys information'' in these statements using tools of information theory (mutual information, relative entropies), but we will not delve on this issue here. The interested reader is referred to any of the many good references on the subject \cite{infotheory}.

From the definition of $\Gamma^{(n)}(P,Q)$, it is clear that generalized digital $n$-circles are combinatorial objects whenever $|P| \geq 2$ or $|Q| \geq 2$, i.e., they represent classes of objects related by combinations and permutations of the symbols in $P$ and $Q$. For example, we can have a $\Gamma^{(3)}(\{1,7\}, \{0,3\})$ = $\langle 1,1,7;\,$ $1,3,1;\,$ $7,1,7 \rangle$ or an equally valid  $\Gamma^{(3)}(\{1,7\}, \{0,3\})$ = $\langle 7,1,7;\,$ $7,0,1;\,$ $1,1,7 \rangle$. This makes generalized digital $n$-circles more ubiquitous than their ``plain'' counterpart.

The generalized Sagan numbers that refer to a generalized digital $n$-circle can be denoted by $S_{\;b}^{(n)}(\alpha;P,Q)$. This number can be defined as indicating the address of the first digit of the first occuring pattern of the class encoded by $\Gamma^{(n)}(P,Q)$. Of course, digital $n$-circles $\Gamma^{(n)}(P,Q)$ can only be found in expansions of irrational numbers in bases $b \geq \max\{\gamma \in P \cup Q\}$, despite the fact that some of the members of the class $\Gamma^{(n)}(P,Q)$ may be found on a given irrational $\alpha$ in base $b$ even when this condition is not obeyed. The reader is invited to parse $\pi$ for the first $\Gamma^{(n)}(P,Q)$ for a few choices of $P$ and $Q$.

\begin{exercise}
How much more frequently should we expect to find $\Gamma^{(n)}(P,Q)$ than the more strict $\Gamma^{(n)}$ in a given $b$-normal number as a funtion of $P$, $Q$, and base $b$?
\end{exercise}

\subsection{BBP-type formul{\ae}}

For some very intersting combination of (somewhat artificially constructed) numbers $\alpha$ and basis $b$, the search for patterns can be done without having to look for the entire expansion of $\alpha$ \cite{crandall,misiu,bailey,bbgal}, since in these cases we know {\it a priori\/} what the sequence of $b$-digits of $\alpha$ looks like, namely,
\[
\alpha = \sum_{k=0}^{\infty}\frac{1}{b^{k}}\frac{p(k)}{q(k)},
\]
with $b$ an integer (the base of the number system) and $p(k)$ and $q(k)$ two polynomials with integer coefficients. With formul{\ae} like this, searches can be done by means of moving windows $(\alpha_{i}\alpha_{i+1}\ldots \alpha_{i+l-1})$ of size $l$ without having to compute all the previous digits up to $\alpha_{i-1}$. This saves time and space in the calculation of the numbers and their digits. Amazingly, $\pi$, $\pi^2$, $\log 2$, $\sqrt{2}$, and a number of related constants can be put in this so called ``BBP-type formula'' in bases 2 and 16 \cite{bailey}. One of the formulas for $\pi$ is truly remarkable and reads
\[
\pi = \sum_{k=0}^{\infty}\frac{1}{16^k}
\left(\frac{4}{8k+1}-\frac{2}{8k+4}-\frac{1}{8k+5}-\frac{1}{8k+6}\right).
\]
It is very tempting, and is a whole research program, to search for similar identities for $\pi$ in other bases. However, J.~Borwein, D.~Borwein, and W.~Galway showed in 2004 that there is no base-$n$ degree-1 BBP-type formula for $\pi$ when $n$ is not a power of two \cite{bbgal}. This result does not rule out the possibility of completely different mathematical approaches that permit one to rapidly calculate digits in other bases starting at an arbitrary point. This is an open problem to which the reader is invited to turn his or her attention!


\section{Conclusion}
\label{conclusion}

Although purely recreational, Sagan numbers hold some properties and involve areas of knowledge that can be used to illustrate several interesting concepts of both pure and applied mathematics, such as those of $b$-ary expansions, normal numbers, basic probability, discrete mathematics, symmetries, and computer graphics, among others.

In order to introduce the numbers, we were led to define a geometrical construction that we called digital $n$-circle. These constructions have a natural habitat in raster graphics, digital imaging processing and pattern recognition. Sagan numbers can also be generalized in many ways by changing the domains of its parameters $\alpha$ and $b$ and also by changing the pattern after which we are looking for, e.g., we may decide to chase for digital $n$-squares instead of digital $n$-circles, or even for the patterns of the ``$\varnothing$'' or the ``$\infty$'' symbols.

If $\pi$ were normal, it could be used as the ultimate encoder of everything symbolic in our lives. For example, we could use its digits to define the $SSN_{16}(\pi)$ as an alternative to our social security number! This is of course a crazy and utterly useless idea, because we would be just changing one code for another, because probably $SSN_{16}(\pi) \gg SSN$ itself, and also because probably it is still unfeasible to find all the necessary digits for everybody within reasonable computational power and time. Moreover, such computations would throw a lot of unnecessary CO$_2$ in the  atmosphere, something we definitely should avoid doing nowadays.

In summary, we hope to have introduced a joyful new class of numbers with many possibilities for the amateur as well as for the professional math\-e\-ma\-ti\-cian---wether of theoretical, applied, or computational inclination---and that one day we could see a table of known Sagan numbers published somewhere just for the wonder of it.

\section*{Acknowledgments}
The author would like to acknowledge F\'{a}bio E. R. Campolim (currently a Ph.D. candidate at UFABC, Brazil) for the many hours of wild conversations (wrestling is a more accurate term, actually) on the ultimate questions and answers as to what life, the Universe and everything else mean.



\begin{thebibliography}{99}

\bibitem{contact}C. Sagan, {\it Contact\/} (New York: Simon \&\ Schuster, 1985).

\bibitem{greene}B. Greene, {\it The elegant Universe: Superstrings, hidden dimensions, and the quest for the ultimate theory\/} (New York: Vintage, 2000).

\bibitem{smolin}L. Smolin, {\it The trouble with Physics: The rise of String Theory, the fall of Science, and what comes next\/} (Boston: Houghton-Mifflin, 2006).

\bibitem{graphics}T. Theoharis, G. Papaioannou, N. Platis, and N. M. Patrikalakis, {\it Graphics and Visualization: Principles \&\ Algorithms\/} (Wellesley: A. K. Peters, 2008).

\bibitem{pitteway}M. L. V. Pitteway, ``Algorithm for drawing ellipses or hyperbol{\ae} with a digital plotter,'' {\it Computer J.\/} {\bf 10}~(3), 282--289 (1967).

\bibitem{vanaken}J. R. Van Aken, ``An efficient ellipse-drawing algorithm,'' {\it IEEE Comput. Graph. \&\ Appl.\/} {\bf 4}~(9), 24--35 (1984); J. Van Aken and M. Novak, ``Curve drawing algorithms for raster displays,'' {\it ACM Trans. Graph.\/} {\bf 4}~(2), 147--169 (1985).

\bibitem{knuth4}D. E. Knuth, {\it The Art of Computer Programming, Volume~4, Fascicle~1: Bitwise Tricks \&\ Techniques, Binary Decision Diagrams\/} (Reading: Addison-Wesley, 2009).

\bibitem{concrete}R. L. Grahan, D. E. Knuth, and O. Patashnik, {\it Concrete Mathematics\/}, 2nd. ed. (Reading: Addison-Wesley, 1998).

\bibitem{lloyd}S. Lloyd, ``Ultimate physical limits to computation,'' {\it Nature\/} {\bf 406} (6799), 1047--1054 (2000); S. Lloyd, ``Computational capacity of the Universe,'' {\it Phys. Rev. Lett.\/} {\bf 88} (23), 237901 (2002); S. Lloyd, {\it Programming the Universe: A Quantum Computer Scientist Takes on the Cosmos\/} (New York: Vintage, 2007).

\bibitem{borel}E. Borel, ``Les probabilit{\'e}s d{\'e}nombrables et leurs applications arithm{\'e}\-ti\-ques,'' {\it Rend. Circ. Mat. Palermo\/} {\bf 27} (1), 247--271 (1909). See also R.~B.~Ash, {\it Real Analysis and Probability\/} (San Diego: Academic Press, 1972), Chapter~7, problem 7.2.7(b).

\bibitem{champer}D. G. Champernowne, ``The construction of decimals normal in the scale of ten,'' {\it J. London Math. Soc.\/} {\bf 8} (4), 254--260 (1933).

\bibitem{erdos}A. H. Copeland and P. Erd\H{o}s, ``Note on normal numbers,'' {\it Bull. Amer. Math. Soc.\/} {\bf 52}, 857--860 (1946).

\bibitem{crandall}D. H. Bailey and R. E. Crandall, ``On the random character of fundamental constant expansions,'' {\it Exp. Math.\/} {\bf 10}~(2), 175--190 (2001).

\bibitem{stone}R. G. Stoneham, ``On $(j,\epsilon)$-normality in the rational fractions,'' {\it Acta Arith.\/} {\bf 16}, 239--237 (1970); R. G. Stoneham, ``A general arithmetic construction of transcendental non-Liouville normal numbers from rational functions,'' {\it Acta Arith.\/} {\bf 16}, 239--253 (1970).

\bibitem{pincus}S. Pincus and B. H. Singer, ``Randomness and degrees of irregularity,'' {\it Proc. Natl. Acad. Sci. USA\/} {\bf 93}~(5), 2083--2088 (1996); S. Pincus and B. H. Singer, ``A recipe for randomness,'' {\it Proc. Natl. Acad. Sci. USA\/} {\bf 95}~(18), 10367--10372 (1998).

\bibitem{kanada}Dr.~Yasumasa Kanada and collaborators announced on Dec.~6, 2002 the record of 1,241,100,000,000 decimal digits calculated for $\pi$ and $1/\pi$. Kanada's laboratory homepage: \url{http://www.super-computing.org}.

\bibitem{daisuke}Dr.~Daisuke Takahashi announced on Aug.~17, 2009 the record of 2,576,980,370,000 decimal digits calculated for $\pi$ and $1/\pi$. Communication at \url{http://www.hpcs.cs.tsukuba.ac.jp/$\sim$daisuke/pi.html}.

\bibitem{wagon}S. Wagon, ``Is Pi normal?,'' {\it Math. Intel.\/} {\bf 7}~(3), 65--67 (1985).

\bibitem{quest}D. Bailey, J. M. Borwein, P. B. Borwein, and S. Plouffe, ``The quest for Pi,'' {\it Math. Intel.\/} {\bf 19}~(1), 50--57 (1997).

\bibitem{calogero}F. Calogero, ``Cool irrational numbers and their rather cool rational approximations,'' {\it Math. Intel.\/} {\bf 25}~(4), 72--76 (2003).

\bibitem{misiu}D. H. Bailey and M. Misiurewicz, ``A strong hot spot theorem,'' {\it Proc. Amer. Math. Soc.\/} {\bf 134}~(9), 2495--2501 (2006).

\bibitem{kapoor}D. H. Bailey, J. M. Borwein, V. Kapoor, and E. N. Weisstein, ``Ten problems in experimental mathematics,'' {\it Am. Math. Month.\/} {\bf 113}~(6), 481--509 (2006).

\bibitem{adams}D. Adams, {\it The Hitchiker's Guide to the Galaxy\/} (UK: Pan Books, 1979); {\it The Restaurant at the End of the Universe\/} (UK: Pan Macmillan, 1980).

\bibitem{infotheory}T. M. Cover and J. A. Thomas, {\it Elements of Information Theory\/}, 2nd.~ed. (Hoboken: John Wiley \&\ Sons, 2006); S. Kullback, {\it Information Theory and Statistics\/} (New York: Dover, 1997).

\bibitem{bailey}D. Bailey, P. Borwein, and S. Plouffe, ``On the rapid computation of various polylogarithmic constants,'' {\it Math. Comput.\/} {\bf 66}~(218), 903--913 (1997).

\bibitem{bbgal}J. M. Borwein, D. Borwein, and W. F. Galway, ``Finding and excluding $b$-ary Machin-type individual digit formulae,'' {\it Can. J. Math.\/} {\bf 56}~(5), 897--925 (2004).

\end{thebibliography}
\end{document}